\documentclass[12pt,draft]{amsart}
\usepackage[all]{xy}
\usepackage{amsfonts}
\usepackage{amssymb}

\usepackage{color}

\newtheorem{teo}{Theorem}[section]
\newtheorem{lem}[teo]{Lemma}

\newtheorem{cor}[teo]{Corollary}

\theoremstyle{definition}
\newtheorem{dfn}[teo]{Definition}
\newtheorem{rk}[teo]{Remark}
\newtheorem{ex}[teo]{Example}

\def\<{\langle}
\def\>{\rangle}
\def\ss{\subset}

\def\r{\rho}

\def\t{\tau}

\def\w{\omega}

\def\G{{\Gamma}}

\def\C{{\mathbb C}}

\def\Z{{\mathbb Z}}

\def\End{\mathop{\rm End}\nolimits}
\def\Ker{\mathop{\rm Ker}\nolimits}
\def\Im{\mathop{\rm Im}\nolimits}
\def\Aut{\operatorname{Aut}}

\def\Coker{\operatorname{Coker}}
\def\Id{\operatorname{Id}}

\def\Tr{\operatorname{Tr}}

\def\1{\mathbf 1}

\def\ola#1{\stackrel{#1}{\longrightarrow}}

\newcommand{\ov}[1]{\overline{#1}}
\newcommand{\til}[1]{\widetilde{#1}}
\newcommand{\wh}[1]{\widehat{#1}}

\def\Fix{\operatorname{Fix}}

\def\N{{\mathbb N}}

\begin{document}

\title[Twisted Burnside-Frobenius theory]
{Twisted Burnside-Frobenius theory\\ for discrete groups}

\author{Alexander Fel'shtyn}
\address{Instytut Matematyki, Uniwersytet Szczecinski,
ul. Wielkopolska 15, 70-451 Szczecin, Poland and Department of Mathematics,
Boise State University,
1910 University Drive, Boise, Idaho, 83725-155, USA}
\email{felshtyn@diamond.boisestate.edu, felshtyn@mpim-bonn.mpg.de}

\author{Evgenij Troitsky}
\thanks{The second author is partially supported by
RFFI Grant  05-01-00923
and Grant ``Universities of Russia''}
\address{Dept. of Mech. and Math., Moscow State University,
119992 GSP-2  Moscow, Russia}
\email{troitsky@mech.math.msu.su}
\urladdr{
http://mech.math.msu.su/\~{}troitsky}

\keywords{Reidemeister number, twisted conjugacy classes,
Burnside-{Frobenius} theorem, solvable group, polycyclic group, conjugacy separable group,
polynomial growth, Osin group}
\subjclass[2000]{20C; 
20E45; 
22D10; 
22D25; 
37C25; 
43A30; 
46L; 
47H10; 
54H25; 
55M20
}

\begin{abstract}
For a wide class of groups including polycyclic and finitely
generated polynomial growth groups it is proved that the Reidemeister
number of an automorphism $\phi$ is equal to the number of
finite-dimensional fixed points of the induced map
$\widehat\phi$ on the unitary dual, if one of these numbers is finite.
This theorem is a natural generalization of the classical
Burnside-Frobenius theorem to
infinite groups. This theorem also has important consequences in
topological dynamics and in some sense is a reply to a remark of J.-P. Serre.
The main technical results proved in the paper yield a tool for a
further progress.
\end{abstract}

\maketitle

\tableofcontents

\section{Introduction}

\begin{dfn}
Let $G$ be a countable discrete group and $\phi: G\rightarrow G$ an
endomorphism.
Two elements $x,x'\in G$ are said to be
 $\phi$-{\em conjugate} or {\em twisted conjugate,}
if and only if there exists $g \in G$ with
$$
x'=g  x   \phi(g^{-1}).
$$
We will write $\{x\}_\phi$ for the $\phi$-{\em conjugacy} or
{\em twisted conjugacy} class
 of the element $x\in G$.
The number of $\phi$-conjugacy classes is called the {\em Reidemeister number}
of an  endomorphism $\phi$ and is  denoted by $R(\phi)$.
If $\phi$ is the identity map then the $\phi$-conjugacy classes are the usual
conjugacy classes in the group $G$.
\end{dfn}

If $G$ is a finite group, then the classical Burnside{-Frobenius}
theorem (see, e.g., \cite{serrerepr},
\cite[p.~140]{Kirillov})
says that the number of
classes of irreducible representations is equal to the number of conjugacy
classes of elements of $G$.  Let $\wh G$ be the {\em unitary dual} of $G$,
i.e. the set of equivalence classes of unitary irreducible
representations of $G$.

If $\phi: G\to G$ is an automorphism, it induces a map $\wh\phi:\wh G\to\wh G$,
$\wh\phi (\r)=\r\circ\phi$.
Therefore, by the Burnside{-Frobenius} theorem, if $\phi$ is the identity automorphism
of any finite group $G$, then we have
 $R(\phi)=\#\Fix(\wh\phi)$.

In \cite{FelHill} it was discovered that
this statement remains true for any automorphism $\phi$ of any finite group $G$.
Indeed, if we
consider an automorphism $\phi$ of a finite group $G$, then $R(\phi)$
is equal to the dimension of the space of twisted invariant functions on
this group. Hence, by Peter-Weyl theorem (which asserts the existence of
a two-side equivariant isomorphism
$C^*(G)\cong \bigoplus_{\r\in\wh G} \End(H_\r)$),
$R(\phi)$ is identified
with the sum of dimensions
$d_\r$ of twisted invariant elements of $\End(H_\r)$, where $\r$ runs over
$\wh G$, and the space of representation $\r$ is denoted by $H_\r$.
By the Schur lemma,
$d_\r=1$, if $\r$ is a fixed point of $\wh\phi$, and is zero otherwise. Hence,
$R(\phi)$ coincides with the number of fixed points of $\wh\phi$.

Dynamical questions have inspired a series of
papers \cite{FelHill,FelBanach,FelshB,FelTro,FelTroVer,FelIndTro},
attempting to generalize this theorem to
the cases of non-identical
automorphisms and of non-finite groups. In these papers
a version of the theorem for almost Abelian groups is proved,
and some examples and particular cases are considered.

In the present paper we introduce the property RP
(Definition \ref{dfn:Rperiodi})
for a countable
discrete group $G$: the $\phi$-class functions of any automorphism $\phi$
with $R(\phi)<\infty$ are periodic in a natural sense.

After some preliminary and technical considerations we prove \textbf{the main
results} of the paper, namely

\begin{enumerate}
    \item {\bf RP respects some extensions:}
Suppose there is an extension $H\to G\to G/H$,
where the group $H$ is a characteristic
{\rm RP}-group; $G/H$ is finitely generated {\rm FC}-group (i.e.
a group with finite conjugacy classes).
Then $G$ is an {\rm RP}-group ({a reformulation of}
Theorem \ref{teo:RPFCexten}).
    \item {\bf Classes of RP groups:} Polycyclic groups and
    finitely generated groups of polynomial growth are RP-groups.
    Moreover, almost-polycyclic groups are RP too.
    (Theorems \ref{teo:polyciclRP}, \ref{teo:polynomRP}, \ref{teo:RPalmostpolyc}).
    The Twisted Burnside-Frobenius
    theorem is valid for them (Theorem \ref{teo:TwistedBurnsPolyn}).
    \item {\bf Twisted Burnside-Frobenius theorem {for RP-groups}:}
Let $G$ be an {\rm RP} group and $\phi$ its automorphism {with $R(\phi)<\infty$}.
Denote by $\wh G_f$ the subset of the unitary dual $\wh G$ related to
finite-di\-men\-si\-on\-al representations. Denote
by $S_f(\phi)$
the number of fixed points of $\wh\phi_f$ on $\wh G_f$. Then
$
R(\phi)=S_f(\phi).
$
(Theorem \ref{teo:twburnsforRP}).
\item { {\bf Twisted Burnside-Frobenius theorem for almost polycyclic groups:}
Let $G$ be an almost polycyclic group. Then $R(\phi)=S_f(\phi)$
if one of these numbers is finite (Theorem \ref{teo:TBFTvdve}).
}
\end{enumerate}

In some sense our theory is a reply to a remark of J.-P.~Serre
\cite[(d), p.34]{serrerepr} that for compact infinite groups,
an analogue of the Burnside-Frobenius theorem is not
interesting: $\infty=\infty$.
It turns out that for infinite discrete groups the situation differs
significantly, and even in non-twisted situations the number of
classes can be finite (for one of the first examples see another
book of J.-P.~Serre \cite{serrtrees}). Several examples
of groups and automorphisms with finite Reidemeister numbers were
obtained and studied in \cite{FelshB,gowon,FelHillWong,FelTroVer,FelIndTro}.

Using the same argument as in \cite{FelTro} one obtains from the
twisted Burnside-Frobenius theorem the following dynamical and number-theoretical
consequence which, together with the twisted Burnside-Frobenius
theorem itself,
is very important for the realization problem of Reidemeister numbers in
topological dynamics and the study of the Reidemeister
zeta-function.

Let $\mu(d)$, $d\in\N$, be the {\em M\"obius function},
i.e.
$$
\mu(d) =
\left\{
\begin{array}{ll}
1 & {\rm if}\ d=1,  \\
(-1)^k & {\rm if}\ d\ {\rm is\ a\ product\ of}\ k\ {\rm distinct\ primes,}\\
0 & {\rm if}\ d\ {\rm is\ not\ square-free.}
\end{array}
\right.
$$

{\sc Congruences for Reidemeister numbers:}
{\it
Let $\phi:G$ $\to G$ be an automorphism of a countable discrete
{\rm RP}-group $G$
such that all numbers $R(\phi^n)$ are finite.
Then one has for all $n$,}
 $$
 \sum_{d\mid n} \mu(d)\cdot R(\phi^{n/d}) \equiv 0 \mod n.
 $$

These theorems were proved previously in a
number of special cases in \cite{FelHill,FelBanach,FelTro,FelTroVer,FelIndTro}.

We would like to emphasize the following important
remarks.

\begin{enumerate}
\item {In the original formulation by Fel'shtyn and Hill \cite{FelHill}
the conjecture about twisted Burnside-Frobenius theorem asserts an equality
of $R(\phi)$ and the number of fixed points of $\wh\phi$ on $\wh G$. This
conjecture was proved in \cite{FelHill,FelTro} for f.g. type I groups.}

\item
In our paper \cite{FelTroVer} with A. Vershik,
we studied a key example which shows that
an RP-group can have infinite-dimensional ``supplementary'' fixed
representations.
Our example was a semi-direct product of the action
of $\Z$ on $\Z\oplus \Z$ by a hyperbolic automorphism. We consider an
 automorphism $\phi$ with finite
Reidemeister number (four to be precise).
$\wh\phi$ has at
least five fixed points on $\wh G$, but exactly four fixed
points on  $\wh G_f$.

This gives a counterexample to the conjecture in its original
formulation (in which we count all fixed points in $\wh G$)
and leads to the formulation, in which we count fixed points only from
$\wh G_f$. This new conjecture is proved in the present paper for
a wide class of f.g. groups.

\item
The extra-fixed-point
phenomenon arises from bad separation properties of $\wh G$ for a general
discrete group $G$. A deeper study leads to the following general theorem.

 {\sc Weak Twisted Burnside theorem} \cite{ncrmkwb}:
{\it
Let $R_*(\phi)$ be the number of Reidemeister classes related to
twisted invariant functions on $G$ from the Fourier-Stieltjes
algebra $B(G)$. Let $S_*(\phi)$ be the number of generalized
fixed points of $\wh\phi$ on the Glimm spectrum of $G$, i.~e., on
the complete regularization of $\wh G$. If one of $R_*(\phi)$ and
$S_*(\phi)$ is finite, then $R_*(\phi)=S_*(\phi)$.
}

The proof is based on a non-commutative version of
the well-known Riesz(-Markov-Kakutani) theorem, which identifies
the space of linear functionals on the algebra $A=C(X)$ with the
space of regular measures on $X$. To prove the Weak Twisted Burnside
theorem we first obtain a generalization of the Riesz theorem to the
case of a non-commutative $C^*$-algebra $A$ using the Dauns-Hofmann
sectional representation theorem (in the same paper \cite{ncrmkwb}).
The corresponding measures
on the Glimm spectrum are functional-valued.
In extreme situations this theorem is tautological. But for
 many cases of group $C^*$-algebras of discrete groups one
obtains some new method for counting twisted conjugacy classes.
This leads to an
approach alternative to the one we present here.

\item
The main Theorem \ref{teo:RPFCexten} allows us to
verify the periodicity of $\phi$-class functions in a number
of cases which are not in the classes described in Section
\ref{sec:RPpolypoly}. Nevertheless for pathological
groups from Section \ref{sec:osin}
even the modified conjecture is not true.
Keeping in mind that for Gromov hyperbolic groups $R(\phi)$ is always infinite
(as well as for
Baumslag-Solitar groups and some generalizations,
cf. \cite{FelPOMI,ll,FelGon,LevittBaums,TabWong})
while in the ``opposite" case the twisted Burnside theorem is proved we
hope that various use of Theorem \ref{teo:RPFCexten} can lead
to a complete resolution of the problem, if the groups from
Section \ref{sec:osin} will be handled.
\end{enumerate}

\medskip
The interest in twisted conjugacy relations has its origins, in particular,
in the Nielsen-Reidemeister fixed point theory (see, e.g. \cite{Jiang,FelshB}),
in Selberg theory (see, eg. \cite{Shokra,Arthur}),
and  Algebraic Geometry (see, e.g. \cite{Groth}).

Concerning some topological applications of our main results, they are already
described in \cite{FelTro}.
The congruences give some necessary conditions for the realization problem
for Reidemeister numbers in topological dynamics.
The relation to Selberg theory will be presented in a forthcoming paper.

\medskip\noindent
{\bf Acknowledgement.} The present research is a
part of our joint research programm
in Max-Planck-Institut f\"ur Mathematik (MPI)
in Bonn.
We would like to thank the MPI for its kind support and
hospitality while the greater part of this work was completed.

The authors are grateful to
R.~Hill,
M.~Sapir,
A.~Shtern,
L.~Vainerman,
A.~Vershik
for helpful discussions and to the referee for valuable comments.

\medskip
The results of Sections \ref{sec:prelim}, \ref{sec:RPpolypoly},
\ref{sec:osin}, and \ref{sec:separat}
are obtained by the authors jointly,
the results of Sections \ref{sec:extens} and
\ref{sec:burnforRP}
are obtained by E.~Troitsky.

\section{Preliminary Considerations}\label{sec:prelim}

The following fact will be useful.
\begin{teo}[{\cite[Theorem 1.41]{Robinson72-1}}]\label{teo:fingeneinind}
If $G$ is a finitely generated group and $H$ is a subgroup with
finite index in $G$, then $H$ is finitely generated.
\end{teo}

\begin{lem}\label{lem:normalcharfiniteind}
Let $G$ be finitely generated, and $H'\ss G$
its subgroup of finite index. Then there is a characteristic
subgroup $H\ss G$ of finite index, $H\ss H'$.
\end{lem}

\begin{proof}
Since $G$ is finitely generated, there is only finitely many
subgroups of the same index as $H'$
(see \cite{Hall}, \cite[\S\ 38]{Kurosh}).
Let $H$ be their intersection. Then $H$ is characteristic, in
particular normal, and of finite index.
\end{proof}

\begin{lem}\label{lem:charabelclassred}
Let $G$ be abelian.
The twisted conjugacy class $H$ of $e$ is a subgroup.
The other classes are cosets $gH$.
\end{lem}

\begin{proof}
The first statement follows from the equalities
$$
h\phi(h^{-1}) g\phi(g^{-1})=gh \phi((gh)^{-1},\quad
(h\phi(h^{-1}))^{-1}=\phi(h) h^{-1}=h^{-1} \phi(h),
$$
where $h\in H$.
For the second statement, suppose $a\sim b$, i.e. $b=ha\phi(h^{-1})$. Then
$$
gb=gha\phi(h^{-1})=h(ga)\phi(h^{-1}), \qquad gb\sim ga.
$$
\end{proof}

Let us denote by $\t_g:G\to G$ the automorphism $\t_g(\til g)=g\til g\,g^{-1}$
for $g\in G$. Its restriction to a normal subgroup we will denote by $\t_g$
as well.

\begin{lem}\label{lem:redklassed}
$\{g\}_\phi k=\{g\,k\}_{\t_{k^{-1}}\circ\phi}$.
\end{lem}

\begin{proof}
Let $g'=f\,g\,\phi(f^{-1})$ be $\phi$-conjugate to $g$. Then
$$
g'\,k=f\,g\,\phi(f^{-1})\,k=f\,g\,k\,k^{-1}\,\phi(f^{-1})\,k
=f\,(g\,k)\,(\t_{k^{-1}}\circ\phi)(f^{-1}).
$$
Conversely, if $g'$ is $(\t_{k^{-1}}\circ\phi)$-conjugate to $g$, then
$$
g'\,k^{-1}=
f\,g\,(\t_{k^{-1}}\circ\phi)(f^{-1})k^{-1}=
f\,g\,k^{-1}\,\phi(f^{-1}).
$$
Hence a shift maps $\phi$-conjugacy classes onto classes related to
another automorphism.
\end{proof}

\begin{cor}\label{cor:rphiequaltaurphi}
$R(\phi)=R(\t_g \circ \phi)$.
\end{cor}

\begin{teo}[see \cite{Jiang}]\label{teo:reidandcokerabel}
Let $A$ be a finitely generated Abelian group, $\psi:A\to A$ its
automorphism. Then $R(\psi)=\#\Coker(\psi-\Id)$, i.e. to the index
of subgroup generated by elements of the form $x^{-1}\psi(x)$.
\end{teo}

\begin{proof}
By Lemma \ref{lem:charabelclassred}, $R(\psi)$ is equal to the index of the subgroup
$H=\{e\}_\psi$. This group consists by definition
of elements of the form $x^{-1}\psi(x)$.
\end{proof}

The following construction relates $\phi$-conjugacy classes and
some conjugacy classes of another group. It was  obtained
in topological context by Boju Jiang and Laixiang
Sun in \cite{Jiang88}.
Consider the action of $\Z$ on $G$, i.e. a homomorphism $\Z\to\Aut(G)$,
$n\mapsto\phi^n$. Let $\G$ be the corresponding semi-direct product
$\G=G\rtimes\Z$:
\begin{equation}\label{eq:defgamma}
    \G:=< G,t\: |\: tgt^{-1}=\phi(g) >
\end{equation}
in terms of generators and relations, where $t$ is a generator of $\Z$.
The group $G$ is a normal subgroup
of $\G$. As a set, $\G$ has the form
\begin{equation}\label{eq:razlgamm}
    \G=\sqcup_{n\in\Z}G\cdot t^n,
\end{equation}
where $G\cdot t^n$ is the coset by $G$ containing $t^n$.

\begin{rk}
Any usual conjugacy class of $\G$ is contained in some $G\cdot t^n$.
Indeed, $g g't^n g^{-1}=gg'\phi^n(g^{-1}) t^n$ and $t g't^n t^{-1}=\phi(g')t^n$.
\end{rk}

\begin{lem}\label{lem:HillRanbij}
Two elements $x,y$ of $G$ are $\phi$-conjugate if and only if $xt$ and $yt$ are conjugate
in the usual sense in $\G$. Therefore $g\mapsto g\cdot t$ is a bijection
from the set of $\phi$-conjugacy classes of $G$ onto the set of conjugacy
classes of $\G$ contained in $G\cdot t$.
\end{lem}

\begin{proof}
If $x$ and $y$ are $\phi$-conjugate then there is a $g\in G$ such that
$gx=y\phi(g)$. This implies $gx=ytgt^{-1}$ and therefore $g(xt)=(yt)g$
so $xt$ and $yt$ are conjugate in the usual sense in $\G$. Conversely,
suppose $xt$ and $yt$ are conjugate in $\G$. Then there is a $gt^n\in\G$
with $gt^nxt=ytgt^n$. From the relation $txt^{-1}=\phi(x)$ we obtain
$g\phi^n(x)t^{n+1}=y\phi(g)t^{n+1}$ and therefore $g\phi^n(x)=y\phi(g)$.
Hence, $y$ and $\phi^n(x)$ are $\phi$-conjugate. Thus,
$y$ and $x$ are $\phi$-conjugate, because $x$ and $\phi(x)$ are always
$\phi$-conjugate: $\phi(x)=x^{-1} x \phi(x)$.
\end{proof}

\section{Extensions and Reidemeister Classes}\label{sec:extens}

Consider a group extension respecting homomorphism $\phi$:
\begin{equation}\label{eq:extens}
 \xymatrix{
0\ar[r]&
H \ar[r]^i \ar[d]_{\phi'}&  G\ar[r]^p \ar[d]^{\phi} & G/H \ar[d]^{\ov{\phi}}
\ar[r]&0\\
0\ar[r]&H\ar[r]^i & G\ar[r]^p &G/H\ar[r]& 0,}
\end{equation}
where $H$ is a normal subgroup of $G$.
The argument below, especially
concerning the role of fixed points, has a partial intersection with
\cite{gowon,goncalves,go:nil1}.

First, notice that the Reidemeister classes of $\phi$ in $G$
are mapped epimorphically onto classes of $\ov\phi$ in $G/H$. Indeed,
\begin{equation}\label{eq:epiofclassforexs}
p(\til g) p(g) \ov\phi(p(\til g^{-1}))= p (\til g g \phi(\til g^{-1}).
\end{equation}
Suppose, $R(\phi)<\infty$. Then the previous remark implies
$R(\ov\phi)<\infty$. Consider a class $D=\{h\}_{\t_g\phi'}$, where
$\t_g(h)=g h g^{-1},$ $g\in G$, $h\in H$. The corresponding equivalence
relation is
\begin{equation}\label{eq:klasstaug}
h\sim \til h h g \phi'(\til h^{-1}) g^{-1}.
\end{equation}
Since $H$ is normal,
the automorphism $\t_g:H\to H$ is well defined.
We will denote
by $D$ the image $iD$ as well. By (\ref{eq:klasstaug})
the shift $Dg$ is a subset of $Hg$, characterized by
\begin{equation}\label{eq:klasstaugsh}
h g\sim \til h (h g) \phi'(\til h^{-1}).
\end{equation}
Hence it is a subset of $\{hg\}_\phi\cap Hg$. The partition
$Hg=\cup (\{h\}_{\t_g \phi' }) g$ is a subpartition of
$Hg=\cup ( Hg\cap \{hg\}_\phi).$

We need the following lemmas.

\begin{lem}\label{lem:ozenkacherezperiod}
Suppose the extension {\rm (\ref{eq:extens})} satisfies
the following conditions:
\begin{enumerate}
    \item $\#\Fix\ov{\phi}=k <\infty$,
    \item $R(\phi)<\infty$.
\end{enumerate}
Then
\begin{equation}\label{eq:ozenreidcherpernonabel}
    R(\phi')\le k\cdot (R(\phi)-R(\ov\phi)+1).
\end{equation}

If $G/H$   is abelian, let $g_i$ be some elements with
$p(g_i)$ being representatives of all
different $\ov\phi$-conjugacy classes, $i=1,\dots,R(\ov\phi)$.
Then
\begin{equation}\label{eq:ozenreidcherperabel}
   \sum_{i=1}^{  R(\ov{\phi})} R(\t_{g_i}\phi')\le k\cdot R(\phi).
\end{equation}
\end{lem}

\begin{proof}
Consider classes $\{z\}_{\phi'}$, $z\in G$, i.e. the classes of
relation $z\sim hz\phi'(h^{-1})$, $h\in H$. The group $G$ acts on
them by $z\mapsto gz\phi(g^{-1})$. Indeed,
\begin{multline*}
    g[\til h h \phi(\til h ^{-1})]\phi(g^{-1})=
(g \til h g^{-1}) (g h \phi(g^{-1}))(\phi(g) \phi(\til h ^{-1})\phi(g^{-1}))\\
    =
(g \til h g^{-1}) (g h \phi(g^{-1})) \phi(g \til hg^{-1})\in
\{g h \phi(g^{-1})\}_{\phi'},
\end{multline*}
because $H$ is normal and $g \til h g^{-1}\in H$. Due to invertibility,
this action of $G$ transposes classes $\{z\}_{\phi'}$ inside one class
$\{g\}_\phi$. Let $d$ denote the number of classes $\{h\}_{\phi'}$
inside $\{h\}_{\phi}\cap H$. Then $d$ does not exceed the number of $g\in G$
such that $p(g)\ov\phi(p(g^{-1}))=\ov e$.
Since two elements $g$ and $gh$
in one $H$-coset induce the same permutation of classes $\{h\}_{\phi'}$,
$d$ does not exceed the number of $z\in G/H$
such that $z\ov\phi(z^{-1})=\ov e$, i.e. $d\le k$. This implies
(\ref{eq:ozenreidcherpernonabel}).

Now we discuss $\phi$-classes over $\ov\phi$-classes other than
$\{\ov e\}_{\ov\phi}$ for an abelian $G/H$. An estimation analogous to
the one above leads to the number of $z\in G/H$
such that $z z_0\ov\phi(z^{-1})=z_0$ for some fixed $z_0$. But for
an Abelian $G/H$ they form the same group $\Fix(\ov\phi)$. This,
together with the description (\ref{eq:klasstaugsh})
of shifts of $D$ at the beginning of this Section, implies
(\ref{eq:ozenreidcherperabel}).
\end{proof}

\begin{lem}\label{lem:fixedofextens}
Suppose, in the extension {\rm (\ref{eq:extens})} the group $H$ is
abelian. Then $\#\Fix(\phi)\le \#\Fix(\phi')\cdot \#\Fix(\ov\phi)$.
\end{lem}

\begin{proof}
Let $s:G/H\to G$ be a section of $p$. If $s(z)h$ is a fixed
point of $\phi$ then
\begin{equation}\label{eq:sotziphiotz}
    (s(z))^{-1}\phi(s(z))=h\phi'(h^{-1}).
\end{equation}
Hence, $z\in \Fix(\ov\phi)$ and left hand side takes $k:=\#\Fix(\ov\phi)$
values $h_1,\dots,h_k$.
Let us estimate the number of $s(z)h$ for a fixed $z$ such that
$(s(z))^{-1}\phi(s(z))=h_i$. These $h$ have to satisfy
(\ref{eq:sotziphiotz}).
Since $H$ is abelian, if one has
$$
h_i=h\phi'(h^{-1})=\til h\phi'(\til h^{-1}),
$$
then $h^{-1}\til h\in \Fix(\phi')$ and we are done.
\end{proof}

\begin{teo}\label{teo:fixandreidforabel}
Let $A$ be a finitely generated Abelian group, $\psi:A\to A$ its
automorphism with $R(\psi)<\infty$. Then $\#\Fix(\psi)<\infty$.

Moreover, $R(\psi)\ge \#Fix(\psi)$.
\end{teo}

\begin{proof}
Let $T$ be the torsion subgroup. It is finite and characteristic.
We obtain the extension $T\to A \to A/T$ respecting $\phi$.
Since $A/T\cong \Z^k$, we have $\Fix(\ov\psi:A/T\to A/T)=\ov e$, by
\cite{Jiang},\cite[Sect. 2.3.1]{FelshB}. Hence, by Lemma
\ref{lem:fixedofextens}, $\#\Fix(\psi)\le \#\Fix(\psi')$, where
$\psi':T\to T$. For any finite abelian group $T$ one
clearly has $\#\Fix(\psi')=R(\psi')$ by Theorem \ref{teo:reidandcokerabel}
(cf. \cite[p.~7]{FelshB}). Finally,
$R(\psi')\le R(\psi)$ by (\ref{eq:ozenreidcherperabel}).
\end{proof}

Recall the following definitions of a class of groups.

\begin{dfn}\label{dfn:FCgroup}
A group with finite conjugacy classes is called FC-\emph{group}.
\end{dfn}

In an FC-group, the elements of finite order form a characteristic
subgroup with locally infinite abelian factor group;
a finitely generated FC-group contains in its center a free abelian
group of finite index in the whole group \cite{BNeumann}.

\begin{lem}\label{lem:fgFChasfinitefix}
An automorphism $\phi$ of a finitely generated {\rm FC}-group $G$
with $R(\phi)<\infty$ has a finite number of fixed points.

The same is true for $\t_x\circ \phi$. Hence, the number of
$g\in G$, such that for some $x\in G$,
$$
g x \phi(g^{-1})=x,
$$
remains finite.
\end{lem}

\begin{proof}
Let $A$ be the center of $G$. As was indicated above, $A$ has finite
index in $G$ and hence, by Theorem \ref{teo:fingeneinind}, is f.g.
Since $A$ is characteristic, one has an extension $A\to G\to G/A$
respecting $\phi$. Hence, $R(\phi')\le R(\phi)\cdot |G/A|$
by Lemma \ref{lem:ozenkacherezperiod}
and $\#\Fix(\phi')\le R(\phi)\cdot |G/A|$ by Theorem
\ref{teo:fixandreidforabel}. Then $\#\Fix(\phi)\le R(\phi)\cdot |G/A|^2$
by Lemma \ref{lem:fixedofextens}.
\end{proof}

\begin{dfn}\label{dfn:Rperiodi}
We say that a group $G$ has the property   {\sc RP} if for
any automorphism $\phi$ with $R(\phi)<\infty$ the
characteristic functions $f$ of {\sc Reidemeister} classes (hence
all $\phi$-central functions) are   {\sc periodic} in the
following sense.

There exists a finite group $K$,
its automorphism $\phi_K$, and epimorphism $F:G\to K$ such that
\begin{enumerate}
    \item The diagram
    $$
    \xymatrix{
G\ar[r]^\phi\ar[d]_F& G\ar[d]^F\\
K\ar[r]^{\phi_K}& K
    }
    $$
    commutes.
    \item $f=F^*f_K$, where $f_K$ is a characteristic function
    of a subset of $K$.
\end{enumerate}

If this property holds for a concrete automorphism $\phi$, we
will denote this property by RP($\phi$).
\end{dfn}

\begin{rk}\label{rk:bijecofclass}
By (2) there is only one class $\{g\}_{\phi}$ which maps onto $\{F(g)\}_{\phi_K}$.
Hence, $F$ induces a bijection of Reidemeister classes.
\end{rk}

\begin{lem}\label{lem:periodihomom}
Suppose, $G$ is f.g.
and $R(\phi)<\infty$. Then
characteristic functions of $\phi$-conjugacy classes
are periodic $($i.e. $G$ has {\rm RP}$(\phi)$ $)$
if and only if their left
shifts generate a finite dimensional space.
\end{lem}

\begin{proof}
By the supposition of finite dimension it follows that the
stabilizer of each $\phi$-conjugacy class has finite index.
Hence, the common stabilizer of all $\phi$-conjugacy classes
under left shifts is an intersection of finitely many subgroups, each
of finite index. Hence, its index is finite.
By Lemma \ref{lem:normalcharfiniteind}
there is some smaller subgroup $G_S$ of finite index
which is normal and $\phi$-invariant.
Then one can take $K=G/G_S$. Indeed, it is sufficient to verify
that the projection $F$ is one to one on classes. In other words,
that each coset of $G_s$ enters only one $\phi$-conjugacy class,
or any two elements of coset are $\phi$-conjugated. Consider
$g$ and $hg$, $g\in G$, $h\in G_s$. Since $h$ by definition preserves
classes, $hg=xg\phi(x^{-1})$ for some $x\in G$, as desired.

Conversely, if  $G$ has {\rm RP}$(\phi)$, the class $\{g\}_\phi$ is
a full pre-image $F^{-1}(S)$ of some class $S\ss K$. Then its left shift
can be described as
\begin{eqnarray*}
g'\{g\}_\phi &=& g'F^{-1}(S)=\{g'g_1 | g_1\in F^{-1}(S)\}=\{g | (g')^{-1}g\in F^{-1}(S)\}\\
&=&
\{g | F((g')^{-1}g)\in S\}=\{g | F(g)\in F(g')(S)\}=F^{-1}(F(g')(S)).
\end{eqnarray*}
Since $K$ is finite, the number of these sets is finite.
\end{proof}

\begin{rk}\label{rk:odnodliavseh}
1) In this situation in accordance with  Lemma \ref{lem:normalcharfiniteind}
the subgroup $G_S$ is characteristic, i.e. invariant under any automorphism.

2) Also, the group $G/G_S$ will serve as $K$
(i.e. give rise to a bijection on sets of Reidemeister classes)
not only for $\phi$ but for
$\t_g\circ \phi$ for any $g\in G$ because they have the same collection
of left shifts of Reidemeister classes by Lemma \ref{lem:redklassed}.
\end{rk}

\begin{teo}\label{teo:RPFCexten}
Suppose, the extension {\rm (\ref{eq:extens})} satisfies
the following conditions:
\begin{enumerate}
    \item $H$ has {\rm RP};
    \item $G/H$ is {\rm FC} f.g.
\end{enumerate}
Then $G$ has {\rm RP}($\phi$).
\end{teo}

\begin{proof}
We have $R(\ov\phi)<\infty$, hence $\#\Fix(\ov\phi)<\infty$
by Lemma \ref{lem:fgFChasfinitefix} as well as $\#\Fix(\t_z\ov\phi)<\infty$ for any $z\in G/H$.
Then by Lemma \ref{lem:ozenkacherezperiod} $R(\t_g\phi')<\infty$ for any $g\in G$.
Let $g_1,...,g_s$, $s=R(\ov\phi)$, be elements of $G$ which are mapped by $p$
to different Reidemeister classes of $\ov\phi$.
Now we
can apply the supposition that $H$ has RP and find a characteristic subgroup
$H_K:=\Ker F\ss H$ of finite index such that $F:H\to K$ gives rise to a bijection
for Reidemeister classes of each of the automorphisms
$\t_{g_i}\circ \phi'$, $i=1,\dots,s$. Moreover,
$H_K$ is contained in the stabilizer (under left shifts)
of each twisted conjugacy class of each of these
automorphisms.
In particular, it is normal in $G$.
Hence, we can take a quotient by $H_K$ of the extension
(\ref{eq:extens}):
$$
\xymatrix{
H\ar@{^{(}->}[r]\ar[d]^F &  G\ar[d]^{F_1} \ar[rd]^p &\\
H/H_K\ar@{^{(}->}[r]\ar@{=}[d]&  G/H_K\ar@{=}[d] \ar[r]^p & G/H\ar@{=}[d]\\
K\ar@{^{(}->}[r]& G_1 \ar[r]^p& \G.
}
$$
The quotient map $F_1: G\to G/H_K$ takes $\{g\}_\phi$ to
$\{g\}_\phi$ and this class is a unique class with this property
(we conserve the notations $e,g,\phi$ for the quotient
objects). Indeed, suppose two classes are mapped onto one.
This means that there are two elements $g$ and $gh_K$ of these {different}
classes,
$h_K\in H_K$. One can choose $\til g\in G$ such that $\til g g \phi(\til g^{-1})=g_i h$
for some $i=1,\dots,s$, $h\in H$. Then
$$
\til g(gh_K)\phi(\til g^{-1}) =\til g \til g^{-1}g_i h \phi(\til g )h_K \phi(\til g^{-1})=
g_i h \wh h_K,\qquad \wh h_K\in H_K.
$$
Hence, $g_i h$ and $g_i  h \wh h_K$ belong to the same (but
different) classes as $g$ and $gh_K$. Moreover, they can not be
$\phi$-conjugate by elements of $H$. Hence (cf.
(\ref{eq:klasstaugsh})), the elements $h$ and $h \wh h_K$ are not
$(\t_{g_i}\circ \phi')$-conjugate in $H$. But $H_K\ni \wh h_K$ is
a subgroup of the intersection of stabilizers (see above). We
arrive to a contradiction.

By Lemma \ref{lem:periodihomom} (applied to $G_1$ and concrete
automorphism $\phi$) for the purpose of finding a map $F_2:G_1\to K_1$
with properties (1) and (2) of the Definition \ref{dfn:Rperiodi}
it is sufficient to verify that shifts of the characteristic
function of $\{h\}_\phi\ss G_1$ form a finite dimensional space, i.e.
the shifts of  $\{h\}_\phi\ss G_1$
form a finite collection of subsets of $G_1$.
After that one can take the composition
$$
G\ola{F_1} G_1 \ola{F_2}  K_1
$$
to complete the proof of theorem.

We were able to apply Lemma \ref{lem:periodihomom} above, because
the group $G_1$ is finitely generated. For example, we can take as
generators all elements of $K$ and some pre-images $s(z_i)\in G_1$
under $p$ of a finite system of generators $z_i$ for $\G$. Indeed,
for any $x\in G_1$ we can find some product of $z_i$ to be equal
to $p(x)$. Then the same product of $s(z_i)$ differs from $x$ by
an element of $K$.

Let us prove that the mentioned space of shifts is finite-dimensional.
By Lemma \ref{lem:redklassed} these shifts of
$\{h\}_\phi\ss G_1$ form a subcollection of
$$
\{x\}_{\t_y\circ \phi},\quad x,y\in G_1.
$$
Hence, by Corollary \ref{cor:rphiequaltaurphi} it is sufficient to
verify that the number of different automorphisms $\t_y:G_1\to G_1$
is finite.

Let $x_1,\dots,x_n$ be some generators of $G_1$. Then the number of different
$\t_y$ does not exceed
$$
\prod_{j=1}^n \# \{\t_y(x_j)\:|\:y\in G_1\}\le
\prod_{j=1}^n |K|\cdot \# \{\t_{z}(p(x_j))\: |\: z\in \G\},
$$
where the last numbers are finite by the definition of FC for $\G$.
\end{proof}

\section{Polycyclic Groups and Groups of Polynomial Growth}\label{sec:RPpolypoly}

Now we will prove using Theorem \ref{teo:RPFCexten}
that some classes of groups are RP groups. As one could expect,
these classes are only
a small part of possible applications of this theorem.

Let $G'=[G,G]$ be the \emph{commutator subgroup} or
\emph{derived group} of $G$, i.e. the subgroup generated by
commutators. $G'$ is invariant under any homomorphism, in
particular it is normal. It is the smallest normal subgroup of $G$
with an abelian factor group. Denoting $G^{(0)}:=G$, $G^{(1)}:=G'$,
$G^{(n)}:=(G^{(n-1)})'$, $n\ge 2$, one obtains \emph{derived series}
of $G$:
\begin{equation}\label{eq:derivedseries}
    G=G^{(0)}\supset G'\supset G^{(2)}\supset \dots\supset G^{(n)}\supset
\dots
\end{equation}
If $G^{(n)}=e$ for some value $n$, i.e. the series (\ref{eq:derivedseries})
    stabilizes by trivial group, then
    the group $G$ is \emph{solvable}.

\begin{dfn}\label{dfn:polycgroup}
A solvable group is a
\emph{polycyclic group}, if it has a derived series with cyclic factors.
\end{dfn}

\begin{teo}\label{teo:polyciclRP}
Any polycyclic group is {\rm RP}.
\end{teo}

\begin{proof}
By Lemma \ref{lem:charabelclassred} any commutative group is RP.
Any extension with $H$ being the commutator subgroup $G'$
of $G$ respects any
automorphism $\phi$ of $G$, because $G'$ is evidently characteristic.
The factor group is abelian, in particular FC.

Any polycyclic group is a result of a finite number of extensions
with finitely generated (cyclic) factor groups,
starting from Abelian group. Thus we can
apply inductively Theorem \ref{teo:RPFCexten}
to complete the proof.
\end{proof}

\begin{teo}\label{teo:nilpotRP}
Any finitely generated nilpotent group is {\rm RP}.
\end{teo}

\begin{proof}
These groups are supersolvable, hence, polycyclic \cite[5.4.6, 5.4.12]{Robinson}.
\end{proof}

\begin{teo}\label{teo:polynomRP}
Any finitely generated group of polynomial growth is {\rm RP}.
\end{teo}

\begin{proof}
By \cite{GromovPolyIHES} a  finitely generated group
of polynomial growth is just
a finite extension of an f.g. nilpotent group $H$.
The subgroup $H$ can be supposed to be characteristic, i.e.
$\phi(H)=H$ for any automorphism $\phi:G\to G$.
Indeed, let $H'\ss G$ be a nilpotent subgroup of index $j$.
Let $H$ be the subgroup from Lemma \ref{lem:normalcharfiniteind}.
By Theorem \ref{teo:fingeneinind} it is finitely generated.
Also, it is nilpotent as
a subgroup of nilpotent group (see \cite[\S\ 26]{Kurosh}).

Since a finite group is a specific case of FC group
and an f.g. nilpotent group has RP by Theorem \ref{teo:nilpotRP},
we can apply Theorem \ref{teo:RPFCexten} to complete the proof.
\end{proof}

\begin{teo}\label{teo:TwistedBurnsPolyn}
The Reidemeister number of any automorphism $\phi$ of an f.g. group
of polynomial growth or polycyclic group {with $R(\phi)<\infty$}
is equal to the
number of finite-dimensional
fixed points of $\wh\phi$ on the unitary dual of this group.
\end{teo}

This theorem will be proved in a more general form below
(Theorem \ref{teo:twburnsforRP}).

\begin{teo}\label{teo:RPalmostpolyc}
Any almost polycyclic group, i.e. an extension of
a polycyclic group with a finite factor group, is {\rm RP}.
\end{teo}

\begin{proof}
The proof repeats almost
literally the proof of Theorem \ref{teo:polynomRP}.
One has to use the fact that a subgroup of a polycyclic group
is polycyclic \cite[p.~147]{Robinson}.
\end{proof}

{These theorems will be proved below in a more strong form
(Theorem \ref{teo:TBFTvdve}).}

\section{The Twisted Burnside-Frobenius Theorem for RP Groups}\label{sec:burnforRP}

\begin{dfn}
Denote by $\wh G_f$ the subset of the unitary dual $\wh G$ related to
finite-di\-men\-si\-on\-al representations.
\end{dfn}

\begin{teo}[Twisted Burnside-Frobenius Theorem for RP-groups]\label{teo:twburnsforRP}
Let $G$ be an {\rm RP} group and $\phi$ its automorphism {with $R(\phi)<\infty$}.
Denote by $S_f(\phi)$
the number of fixed points of $\wh\phi_f$ on $\wh G_f$. Then
$$
R(\phi)=S_f(\phi).
$$
\end{teo}

\begin{proof}
The coefficients of finite-dimensional non-equivalent irreducible
representations of $G$ are linear independent by Frobenius-Schur theorem
(see \cite[(27.13)]{CurtisReiner}). Moreover, the coefficients of non-equivalent unitary
finite-dimensional irreducible representations are orthogonal to each other
as functions on the universal compact
group associated with the initial group \cite[16.1.3]{DixmierEng} by the
Peter-Weyl theorem. Hence, their linear combinations are orthogonal to each other
as well.

It is sufficient to verify the following three statements:

1) If $R(\phi)<\infty$, then each $\phi$-class function is a finite
linear combination of twisted-invariant functionals being coefficients
of points of $\Fix\wh\phi_f$.

2) If $\rho\in\Fix\wh\phi_f$, there exists one and only one
(up to scaling)
twisted invariant functional on $\rho(C^*(G))$ (this is a finite
full matrix algebra).

3) For different $\rho$ the corresponding $\phi$-class functions are
linearly independent. This follows from the remark at the beginning of
the proof.

Let us remark that the property RP implies in particular that
$\phi$-central functions (for $\phi$ with $R(\phi)<\infty$) are
functionals on $C^*(G)$, not only $L^1(G)$, i.e. are in the
Fourier-Stieltijes algebra $B(G)$.

The statement 1) follows from the RP property. Indeed,
this $\phi$-class function $f$ is a linear combination of functionals
coming from some finite collection $\{\r_i\}$ of elements
of $\wh G_f$ (these
representations $\r_1,\dots,\r_s$ are in fact representations of
the form $\pi_i\circ F$, where $\pi_i$ are irreducible representations
of the finite group $K$ and $F:G\to K$, as in the definition of RP).
So,
$$
f=\sum_{i=1}^s f_i\circ \r_i,\quad \r_i:G\to \End (V_i),
\quad f_i:\End (V_i)\to \C, \quad \r_i\ne \r_j,\: (i\ne j).
$$
For any $g, \til g \in G$ one has
$$
\sum_{i=1}^s f_i (\r_i(\til g ))=f(\til g )=
f(g\til g \phi(g^{-1}))=\sum_{i=1}^s f_i (\r_i(g\til g \phi(g^{-1}))).
$$
By the observation at the beginning of the proof concerning
linear independence,
$$
f_i (\r_i(\til g ))=f_i (\r_i(g\til g \phi(g^{-1}))).\qquad i=1,\dots,s,
$$
i.e. $f_i$ are twisted-invariant.  For any $\r\in \wh G_f$,
$\r:G\to \End (V)$, any
functional $\w: \End (V)\to \C$
has the form $a\mapsto \Tr(ba)$ for some fixed $b\in \End (V)$.
Twisted invariance implies twisted invariance of $b$ (evident details can
be found in \cite[Sect. 3]{FelTro}). Hence, $b$ is intertwining between
$\r$ and $\r\circ\phi$ and $\r\in \Fix(\wh\phi_f)$. The uniqueness of
intertwining operator (up to scaling) implies 2).
\end{proof}

\section{Counterexamples}\label{sec:osin}
Now let us consider some counterexamples to this statement
for pathological (monster) discrete groups.
Suppose, an infinite discrete group
$G$ has a finite number of conjugacy classes.
Such examples can be found in \cite{serrtrees} (HNN-group),
\cite[p.~471]{olsh} (Ivanov group), and \cite{Osin} (Osin group).
Then evidently, the characteristic function of unity element is not
almost-periodic and the argument above is not valid. Moreover, let us
show, that these groups produce counterexamples to the above theorem.

\begin{ex}\label{ex:osingroup}
For the Osin group the Reidemeister number $R(\Id)=2$,
while it has only one finite-dimensional irreducible
representation (trivial 1-dimensional).
Indeed, the Osin group is
an infinite, finitely generated group $G$ with exactly two conjugacy classes.
All nontrivial elements of this group $G$ are conjugate. So, the group $G$
is simple, i.e. $G$ has no nontrivial normal subgroup.
This implies that group $G$ is not residually finite
(by definition of residually finite group). Hence,
it is not linear (by Mal'cev theorem \cite{malcev}, \cite[15.1.6]{Robinson})
and has no finite-dimensional irreducible unitary
representations with trivial kernel. Hence, by simplicity of $G$, it has no
finite-dimensional
irreducible unitary representation with nontrivial kernel, except of the
trivial one.

Let us remark that the Osin group is non-amenable, contains the free
group in two generators $F_2$,
and has exponential growth.
\end{ex}

\begin{ex}\label{ex;ivanovgroup}
For large enough prime numbers $p$,
the first examples of finitely generated infinite periodic groups
with exactly $p$ conjugacy classes were constructed
by Ivanov as limits of hyperbolic groups (although hyperbolicity was not
used explicitly) (see \cite[Theorem 41.2]{olsh}).
The Ivanov group $G$ is an infinite periodic
2-generator  group, in contrast to the Osin group, which is torsion free.
The Ivanov group $G$ is also a simple group.
The proof (kindly explained to us by M. Sapir) is the following.
Denote by $a$ and $b$ the generators of $G$ described in
\cite[Theorem 41.2]{olsh}.
In the proof of Theorem 41.2 on  \cite{olsh}
it was shown that each of elements of $G$
is conjugate in $G$ to a power of generator $a$ of order $s$.
Let us consider any normal subgroup $N$ of $G$.
Suppose
$\gamma \in N$. Then $\gamma=g a^sg^{-1}$ for some $g\in G$ and some $s$.
Hence,
$a^s=g^{-1} \gamma g \in N$ and from  periodicity of $a$, it follows that also
$ a\in N$
as well as $ a^k \in N$  for any $k$, because $p$ is prime.
Then any element $h$ of $G$ also belongs to $N$
being of the form $h=\til h a^k (\til h)^{-1}$, for  some $k$, i.e., $N=G$.
Thus, the group $G$ is simple. The discussion can be completed
in the same way as in the case of the Osin group.
\end{ex}

\begin{ex}
In paper \cite{HNN}, Theorem III and its corollary,
G.~Higman, B.~H.~Neumann, and H.~Neumann
proved that any locally infinite countable group $G$
can be embedded into a countable group $G^*$ in which all
elements except the unit element are conjugate to each other
(see also \cite{serrtrees}).
The discussion above related Osin group remains valid for $G^*$
groups.
\end{ex}

Let us remark that almost polycyclic groups are residually finite
(see e.g. \cite[5.4.17]{Robinson})
while the groups from these counterexamples are not residually finite,
as it was shown by definition. Hence, we would like to
complete this section with the following question.

\noindent{\bf Question.}
Suppose $G$ is a residually finite group and $\phi$ is its endomorphism
with finite $R(\phi)$. Does $R(\phi)$ equal $S_f(\phi)$?

\section{Twisted Conjugacy Separateness}\label{sec:separat}
In fact, the notion of RP group is closely related to a generalization
of the following well-known notion.

\begin{dfn}
A group $G$ is \emph{conjugacy separable} if any pair $g$, $h$ of
   non-conjugate elements of $G$ are non-conjugate in some finite
   quotient of $G$.
\end{dfn}

It was proved that almost polycyclic groups
are conjugacy separable (\cite{Remes,Form}, see also \cite[Ch.~4]{DSegalPoly}).
Also, residually finite recursively presented
Burnside $p$-groups constructed by R. I. Grigorchuk
\cite{GrFA} and by N. Gupta and S. Sidki \cite{GuSi}
are shown to be conjugacy separable
when $p$ is an odd prime in \cite{WiZa}.

We can introduce the following notion, which coincides with the previous
definition in the case $\phi=\Id$.

\begin{dfn}
A group $G$ is \emph{$\phi$-conjugacy separable} with respect to
an automorphism $\phi:G\to G$ if any pair $g$, $h$ of
   non-$\phi$-conjugate elements of $G$ are non-$\ov\phi$-conjugate in some finite
   quotient of $G$ respecting $\phi$.
\end{dfn}

One gets immediately the following statement.

\begin{teo}\label{teo:phiconjRP}
Suppose, $R(\phi)<\infty$. Then $G$ is $\phi$-conjugacy separable if and
only if $G$ is {\rm RP}$(\phi)$.
\end{teo}

\begin{proof}
Indeed, let $F_{ij}:G\to K_{ij}$ distinguish $i$th and $j$th $\phi$-conjugacy
classes, where $K_{ij}$ are finite groups, $i,j=1,\dots,R(\phi)$. Let
$F:G\to \oplus_{i,j}K_{ij}$, $F(g)=\sum_{i,j}F_{ij}(g)$, be the diagonal
mapping and $K$ its image. Then $F:G\to K$ gives RP$(\phi)$.

The opposite implication is evident.
\end{proof}

\begin{teo}\label{teo:razdelgamg}
Let $F:\G\to K$ be a morphism onto a finite group $K$ which separates
two conjugacy classes of $\G$ in $G\cdot t$. Then the restriction
$F_G:=F|_G:G\to \Im(F|_G)$
separates the corresponding $($by the bijection from
Lemma {\rm \ref{lem:HillRanbij})} $\phi$-conjugacy classes in $G$.
\end{teo}

\begin{proof}
First let us remark that $\Ker(F_G)$ is $\phi$-invariant.
Indeed, suppose $F_G(g)=F(g)=e$. Then
$$
F_G(\phi(g))=F(\phi(g))=F(tgt^{-1})=F(t)F(t)^{-1}=e
$$
(the kernel of $F$ is a normal subgroup).

Let $gt$ and $\til gt$ be some representatives of the mentioned
conjugacy classes. Then
$$
F((ht^n)gt (ht^n)^{-1})\ne F(\til gt),\qquad\forall h\in G,\: n\in \Z,
$$
$$
F(ht^ngt )\ne F(\til gt ht^n),\qquad\forall h\in G,\: n\in \Z,
$$
$$
F(h\phi^n(g)t^{n+1} )\ne F(\til g \phi(h)t^{n+1}),\qquad\forall h\in G,\: n\in \Z,
$$
$$
F(h\phi^n(g))\ne F(\til g \phi(h)),\qquad\forall h\in G,\: n\in \Z,
$$
in particular, $F(hg\phi(h^{-1}))\ne F(\til g )$ $\forall h\in G$.
\end{proof}

\begin{teo}\label{teo:conjsepandRP}
Let some class of conjugacy separable groups be closed under
taking semidirect products by $\Z$. Then this class consists of
{\rm RP} groups.
\end{teo}

\begin{proof}
This follows immediately from Theorem \ref{teo:razdelgamg} and
Theorem \ref{teo:phiconjRP}.
\end{proof}

As an application we obtain another proof of the main theorem for
almost polycyclic groups.

\begin{teo}
Any almost polycyclic group is an {\rm RP} group.
\end{teo}

\begin{proof}
The class of almost polycyclic groups is closed under
taking semidirect products by $\Z$. Indeed, let $G$ be an
 almost polycyclic group. Then there exists a characteristic
 (polycyclic) subgroup $P$ of finite index in $G$. Hence,
 $P\rtimes \Z$ is a polycyclic normal group of $G\rtimes \Z$
 of the same finite index.

Almost polycyclic groups are
conjugacy separable (\cite{Remes,Form}, see also \cite[Ch.~4]{DSegalPoly}).
It remains to apply Theorem \ref{teo:conjsepandRP}.
\end{proof}

\begin{rk}
It is clear that this approach (i.e. the using of Theorem \ref{teo:conjsepandRP})
can be extended to a number of situations (cf. \cite{RiSegZal98}).
The results of the present paper are generalized to
the case of two endomorphisms and their coincidences in \cite{bitwisted}.
The method of proof there is that of Sections \ref{sec:extens} -- \ref{sec:burnforRP}.
The method of Section \ref{sec:separat} does not work in this case.

In fact, in \cite{bitwisted} it is obtained in particular a generalization of
the main result to the case of \emph{endomorphism.} For the method of Sections
\ref{sec:extens} -- \ref{sec:burnforRP} the condition of finiteness of Reidemeister
number is very important. On the other hand for the method of Section \ref{sec:separat}
the property of $\phi$ to be an automorphism is very important. Moreover, we can obtain
the following final form of the twisted Burnside-Frobenius theorem for almost
polycyclic and some other classes of groups.
\end{rk}

\begin{teo}[Twisted Burnside-Frobenius theorem for $\phi$-conjugacy separable groups]\label{teo:TBFTvdve}
Let $G$ be an almost polycyclic group $($or, more generally, any group of a
class under the hypothesis of Theorem {\rm\ref{teo:conjsepandRP}}, or even more
generally, a $\phi$-conjugacy separable group$)$. Then $R(\phi)=S_f(\phi)$
if one of these numbers is finite.
\end{teo}

\begin{proof}
It remains to prove that $R(\phi)<\infty$ if $S_f(\phi)<\infty$.
By the definition of a $\phi$-conjugacy separable group
the Reidemeister classes of $\phi$ can be separated by maps to finite groups.
Hence, taking representations of these finite groups and applying
the twisted Burnside-Frobenius theorem to these groups we obtain that
for any pair of Reidemeister classes there exists a function being a
coefficient of a finite-dimensional unitary representation, which
distinguish these classes. Hence,
if $R(\phi)=\infty$, then there are infinitely
many linearly independent twisted invariant functions being
coefficients of finite dimensional
representations. But there are as many such functionals, as $S_f(\phi)$.
\end{proof}
\begin{rk}
As it is clear from the consideration in Section \ref{sec:osin}, the Osin group
and its predecessors are not conjugacy separable. In this relation it is
natural to reformulate the question from Section \ref{sec:osin} in the following
form.
\end{rk}

\noindent{\bf Question.}
Suppose $G$ is a residually finite group and $\phi$ is its endomorphism
with finite $R(\phi)$. Does this imply that $G$ is $\phi$-conjugacy separable ?

Some other aspects of the twisted analog of Burnside-Frobenius theory are studied in
\cite{frobfin}.


\def\cprime{$'$} \def\dbar{\leavevmode\hbox to 0pt{\hskip.2ex \accent"16\hss}d}
  \def\polhk#1{\setbox0=\hbox{#1}{\ooalign{\hidewidth
  \lower1.5ex\hbox{`}\hidewidth\crcr\unhbox0}}}
\providecommand{\bysame}{\leavevmode\hbox to3em{\hrulefill}\thinspace}
\providecommand{\MR}{\relax\ifhmode\unskip\space\fi MR }
\providecommand{\MRhref}[2]{%
  \href{http://www.ams.org/mathscinet-getitem?mr=#1}{#2}
}
\providecommand{\href}[2]{#2}

\end{document}